\documentstyle[11pt]{article}

\newfam\frakfam
\font\teneufm=eufm10     \textfont\frakfam=\teneufm
\font\seveneufm=eufm7    \scriptfont\frakfam=\seveneufm
\font\fiveeufm=eufm5     \scriptscriptfont\frakfam=\fiveeufm
\def\frak{\fam\frakfam \teneufm}

\newfam\msbfam
\font\tenmsb=msbm10      \textfont\msbfam=\tenmsb
\font\sevenmsb=msbm7     \scriptfont\msbfam=\sevenmsb
\font\fivemsb=msbm5      \scriptscriptfont\msbfam=\fivemsb
\def\Bbb{\fam\msbfam \tenmsb}

\font\titlebf=cmbx10  scaled \magstep2

\newcommand{\ie}{\mbox{\em i.e.}}
\newcommand{\fq}{\mbox{${\Bbb F}_q\,$}}
\newcommand{\fp}{\mbox{${\Bbb F}_p\,$}}
\newcommand{\clo}{\mbox{${\overline {\Bbb F}}_p\,$}}
\newcommand{\nq}{\mbox{${\Bbb F}_{q^n}\,$}}

\newcommand{\sn}{\mbox{${\frak S}_n\,$}}
\newcommand{\op}{\mbox{${\cal O}_P\,$}}

\newcommand{\pf}{\mbox{${\Bbb P}(F)\,$}}

\newcommand{\pti}{\mbox{${\Bbb P}(T_i)\,$}}
\newcommand{\jth}{\mbox{$j^{th}\;$}}
\newcommand{\ith}{\mbox{$i^{th}\;$}}

\newcommand{\degdiff}{\mbox{\rm degDiff}}
\newcommand{\gal}{\mbox{\rm Gal}}

\newtheorem{theorem}{Theorem}[section]
\newtheorem{lemma}[theorem]{Lemma}
\newtheorem{proposition}[theorem]{Proposition}
\newtheorem{corollary}[theorem]{Corollary}
\newtheorem{definition}[theorem]{Definition}

\newtheorem{example}[theorem]{\sc Example}

\begin{document}
\font\titlebf=cmbx10 scaled \magstep2
\centerline{\titlebf Towers of function fields with extremal properties}
\medskip
\centerline{{\bf Vinay Deolalikar} \footnote[1]{ This work forms part of the author's doctoral research and was supervised by the late Prof. Dennis Estes.}}
\medskip
\medskip
\centerline{\sc Dedicated to the late Prof. Dennis Estes}

\section{Introduction}

For $F/K$ an algebraic function field in one variable over a finite field of constants $K$ (\ie, $F$ is a finite algebraic extension of $K(x)$ where $x \in F$ is transcendental over $K$), let $N(F)$ and $g(F)$ denote the number of places of degree one and the genus, respectively, of $F$.

Let ${\cal F} = (F_1,F_2,F_3,\ldots)$ be a tower of function fields, each defined over $K$. Further, we will assume that $F_1 \subseteq F_2 \subseteq F_3 \ldots$, where $F_{i+1}/F_i$ is a finite separable extension and  $g(F_i) > 1$ for some $i \geq 1$. This follows the conventions of \cite{GarSti2}.

In this paper, the techniques developed in \cite{DeoEst1} and \cite{DeoEst2} are applied to splitting rational places in towers of function fields. While the basic ideas are the same, it has to be kept in mind that what is optimal at one stage of the tower may lead to complications at later stages.

Let ${\cal F}$ be as above. It is known that the sequence $(N(F_i)/g(F_i))$ converges as $i \rightarrow \infty$ \cite{GarSti2}. Let $\lambda({\cal F}) := \lim_{i \rightarrow \infty} N(F_i)/g(F_i)$.

There are known bounds on the behaviour of function fields over a finite field \fq. Let $N_q(g)$ := max\{$N(F) | F$ a function field over ${\Bbb F}_q$ of genus $g(F) = g$\}. Also, let
\begin{equation}
A(q) := \limsup_{g \rightarrow \infty} N_q(g)/g,
\end{equation}
then the Drinfeld-Vladut bound \cite{DriVla1} says that
\begin{equation}
A(q) \leq \sqrt{q} - 1.
\end{equation}
Ihara \cite{Iha1}, and Tsafasman, Vladut and Zink \cite{TsfVlaZin1} showed that this bound can be met in the case where $q$ is a square. It is not known what the value of $A(q)$ is for non-square $q$, though there are results by Serre \cite{Ser1,Ser2,Ser3} and Schoof \cite{Sch1} in this direction.

Clearly, for a tower of function fields ${\cal F} = (F_1,F_2,\ldots)$, $F_i/\fq$, we have that
\begin{equation}
 0 \leq \lambda({\cal F}) \leq A(q).
\end{equation}

Garcia and  Stichtenoth \cite{GarSti1, GarSti2} gave two explicitly constructed towers of function fields over a field of square cardinality that meet the Drinfeld-Vladut bound. In \cite{GarSti3}, they gave more explicit descriptions of towers of function fields over \fq, with $\lambda({\cal F}) > 0$. These also meet the Drinfeld-Vladut bound in some cases where the underlying field of constants is of square cardinality.

Elkies, in \cite{Elk1}, gave eight explicit iterated equations for towers of modular curves, which also attained the Drinfeld-Vladut bound over certain fields and showed that the examples presented in \cite{GarSti1} and \cite{GarSti3} were also modular. He then conjectured that all asymptotically optimal towers would, similarly, be modular.

In \cite{DeoEst1}, the author used the notion of symmetry of functions to describe explicitly constructed extensions of function field in which all rational places except one split completely. In \cite{DeoEst2}, it was shown that on generalizing the notion of symmetry to include the so-called ``quasi-symmetric'' functions, one could actually split all the rational places in an extension of function fields. Furthermore, in both these cases, infinite families of extensions with such properties were obtained.

In this paper, techniques developed in \cite{DeoEst1} and \cite{DeoEst2} are applied to the problem of splitting rational places in a tower of function fields. Towards that end, infinite families of towers in which all the rational places split completely throughout the tower are described. Infinite families of towers in which all rational places, except one, split completely throughout the tower are also described. It is observed that inspite of such splitting behaviour at the rational places, all these towers have $\lambda({\cal F}) = 0$. In that sense, the main accent here is not so much on obtaining a high value for  $\lambda({\cal F})$, as it is to show the existence of certain explicitly constructed families of towers in which all rational places split completely throughout the tower. In addition, it is hoped that these examples will lead to a better general understanding of what makes  $\lambda({\cal F}) > 0$.  Two examples of towers with  $\lambda({\cal F}) > 0$ presented in \cite{GarSti3} are also generalized, resulting in  infinite families of such towers. Subfamilies of these attain the Drinfeld-Vladut bound.

\section{Notation}

For symmetric polynomials:
\begin{tabbing}
\sn \hspace{1.6cm} \= the symmetric group on $n$ characters \\
${\bf s}_{n,i}(X)$ \> the \ith elementary symmetric polynomial on $n$ variables \\
$q$  \>  a power of a prime $p$ \\
${\Bbb F}_l$  \> the finite field of cardinality $l$ \\
$s_{n,i}(t)$ \> the \ith $(n,q)$-elementary symmetric polynomial
\end{tabbing}

For function fields and their symmetric subfields:
\begin{tabbing}
$K$ \hspace{1.6cm} \= the finite field of cardinality $q^n$, where $n>1$ \\
$F/K$ \> an algebraic function field in one variable whose full field of constants is $K$ \\
$F_s$ \> the subfield of $F$ comprising $(n,q)$-symmetric functions \\
$F_s^\phi$ \> the subfield of $F_s$ comprising functions whose coefficients are from \fq \\
$F_{qs}$ \> the subfield of $F$ comprising $(n,q)$-quasi-symmetric functions \\
$F^{\phi}_{qs}$ \> the subfield of $F_{qs}$ comprising functions whose coefficients are from \fq \\
$E$ \> a finite separable extension of $F$, $E=F(y)$ where $\varphi(y) = 0$
for some irreducible  \\
 \> polynomial $\varphi[T] \in F[T]$
\end{tabbing}

For a generic function field $F$:
\begin{tabbing}
${\Bbb P}(F)$  \hspace{1.3cm} \=  the set of places of $F$ \\
$N(F)$  \> the number of places of degree one in $F$ \\
$g(K)$ \> the genus of $F$   \\
$P$ \> a generic place in $F$  \\
$v_P$ \> the normalized discrete valuation associated with the place $P$ \\
\op \> the valuation ring of the place $P$ \\
$P'$ \> a generic place lying above $P$ in a finite separable extension of $F$\\
$e(P'|P)$ \> the ramification index for $P'$ over $P$ \\
$d(P'|P)$ \> the different exponent for $P'$ over $P$
\end{tabbing}

For the rational function field $K(x)$:
\begin{tabbing}
$P_\alpha$ \hspace{1.6cm} \= the place in $K(x)$ that is the unique zero of $x-\alpha, \; \alpha \in K$ \\
$P_\infty$ \> the place in $K(x)$ that is the unique pole of $x$
\end{tabbing}

For towers of function fields:
\begin{tabbing}
${\cal F}$ \hspace{1.8cm} \= a tower of function fields $F_i \subseteq F_2 \subseteq F_3 \ldots$ \\
$\lambda(\cal F)$ \> $\displaystyle\lim_{i \rightarrow \infty}(N(F_i)/g(F_i))$
\end{tabbing}

\section{Preliminaries}

In this section we state some preliminary results. For detailed proofs of these,  please refer \cite{DeoEst1} and \cite{DeoEst2}.

\begin{proposition} \label{proposition:Art-Sch}
Let $F/K$ be an algebraic function field, where $K=\nq$ is algebraically closed in $F$. Let $w \in F$ and assume that there exists a place $P \in \pf$ such that
$$ v_P(w) = -m, \, m > 0 \mbox{  and  } \gcd(m,q)  = 1.  $$
Then the polynomial $l(T)-w = a_{n-1}T^{q^{n-1}} + a_{n-2}T^{q^{n-2}}+ \ldots + a_0T - w \in F[T] $ is absolutely irreducible. Further, let $l(T)$ split into linear factors over $K$. Let $E=F(y)$ with
$$  a_{n-1}y^{q^{n-1}} + a_{n-2}y^{q^{n-2}}+ \ldots + a_0y = w. $$
Then the following hold:
\begin{enumerate}
\item $E/F$ is a Galois extension, with degree $[E:F] = q^{n-1}$. $\gal(E/F) = \{\sigma_\beta:y \rightarrow y + \beta\}_{l(\beta)=0}$.
\item $K$ is algebraically closed in $E$.
\item The place $P$ is totally ramified in $E$. Let the unique place of $E$ that lies above $P$ be $P'$. Then the different exponent $d(P'|P)$ in the extension $E/F$ is given by
$$ d(P'|P) = (q^{n-1}-1)(m+1).$$
\item Let $R \in \pf$, and $v_R(w) \geq 0$. Then $R$ is unramified in $E$.
\item If $a_{n-1}=\ldots=a_0=1$, and if $Q\in \pf$ is a zero of $w-\gamma$, with $\gamma \in \fq$. Then $Q$ splits completely in $E$.
\end{enumerate}
\end{proposition}
{\bf Proof}. For (i) - (iv), pl. refer \cite{Sul1}. For (v), notice that under the hypotheses, the equation $ T^{q^{n-1}} +T^{q^{n-2}}+ \ldots + T = \gamma $ has $q^{n-1}$ distinct roots in $K$.  \hfill $\Box$

For many of the extensions that we will describe, there exists no place where the hypothesis of Proposition~\ref{proposition:Art-Sch} is satisfied, namely, that the valuation of $w$ at the place is negative and coprime to the characteristic. In particular, we need a criterion for determining the irreducibility of the equations that we will need to use. We provide such a criterion in Proposition~\ref{proposition:irreducibility} and Corollary~\ref{corollary:irreducibility}.

\begin{proposition} \label{proposition:irreducibility}
Let $V$ be a finite subgroup of the additive group of \clo. Then $V$ is a \fp-vector space. Define $L_V(T) = \prod_{v \in V}(T-v)$. Thus, $L_V(T)$ is a separable \fp-linear polynomial whose degree is the cardinality of $V$. Now let $h(T,x) =  L_V(T) - f(x)$, where $f(x) \in \clo[x]$. Then, $h(T,x)=L_V(T) - f(x)$ is reducible over $\clo[T,x]$ iff there exists a polynomial $g(x) \in \clo[x]$ and a proper additive subgroup $W$ of $V$ such that $f(x) = L_{W'}(g(x))$, where $W' = L_W(V)$.
\end{proposition}
For a proof of this proposition, please refer to \cite{DeoThe} or \cite{DeoEst1}.

\begin{definition}
For $f(x) \in \clo[x]$, a coprime term of $f$ is a term with non-zero coefficient in $f$ whose degree is coprime to $p$. The coprime degree of $f$ is the degree of the coprime term of $f$ having the largest degree.
\end{definition}

\begin{corollary} \label{corollary:irreducibility}
Let $f(x) \in \clo[x]$. Let there be a coprime term in $f(x)$ of degree $d$, such that there are no terms of degree $dp^i$ for $i>0$ in $f(x)$. Then $L_V(T) -f(x)$ is irreducible for any subgroup $V \subset \clo$.
\end{corollary}
{\bf Proof}. Suppose $f(x)$ is the image of a linear polynomial $\sum a_nx^{p_n}$. Then the coprime term can only occur in the image of the term $a_0x$. But then, the images of the coprime term under $a_nx^{p^n}$, for $n>0$ will have degrees that contradict the hypothesis.

\begin{lemma} \label{lemma:subextensions}
 Let $F = K(x)$, where $K=\nq, q = p^m, r=m(n-1),$ and $E = F(y)$, where $y$ satisifes the following equation:
$$ y^{q^{n-1}} + y^{q^{n-2}} + \ldots + y = f(x), $$
and $f(x) \in F$ is not the image of any element in $F$ under a linear polynomial.
Then the following hold:
\begin{enumerate}
\item $E/F$ is a Galois extension of degree $[E:F]=q^{n-1}$. $\gal(E/F) = \{\sigma_\beta: y \rightarrow y + \beta \}_{s_{n,1}(\beta)=0}$ can be identified with the set of elements in \nq whose trace in \fq is zero by $\sigma_\beta \leftrightarrow \beta$. This gives it the structure of a $r$-dimensional \fp vector space.
\item There exists a (non-unique) tower of subextensions
$$ F=E^0 \subset E^1 \subset \ldots \subset E^{r} = E, $$
such that for $0 \leq i \leq r-1,\; [E^{i+1}: E^i]$ is a Galois extension of degree $p$.
\item Let $\{b_i\}_{1 \leq i \leq r}$ be a \fp-basis for $\gal(E/F)$. Then we can build one tower of subextensions as in {\rm (ii)} as follows. We set $E^j$ to be  the fixed field of the subgroup of the Galois group that corresponds to the \fp-subspace generated by $\{b_1,b_2,\ldots, b_{r-j}\}$. Then, the generators of $E^{j}$ are $\{y_1,y_2,\ldots,y_j\}$, where $y_1,y_2,\ldots,y_{r}=y$ satisfy the following relations:
\begin{eqnarray*}
       y^p - B_{r}^{p-1} y &=& y_{r-1}, \\
       y_{r-1}^p - B_{r-1}^{p-1} y_{r-1} &=& y_{r-2}, \\
          \vdots &  & \vdots \\
       y_1^p - B_{1}^{p-1}y_1 &=& f(x),
\end{eqnarray*}
where,
$$
\begin{array}{rcll}
  \beta_{r,j} &=& b_{r-j+1},  & B_{r} = \beta_{r,r}, \\
  \beta_{r-1,j} &=& \beta_{r,j}^p - B_{r}^{p-1}\beta_{r,j}, & B_{r-1} = \beta_{r-1,r-1},\\
 \vdots & & \vdots &\vdots  \\
  \beta_{1,j}   &=& \beta_{2,j}^p - B_2^{p-1}\beta_{2,j}, &  B_1 = \beta_{1,1}.
\end{array}
$$
\end{enumerate}
\end{lemma}
For a proof of this lemma, please refer to \cite{DeoThe} or \cite{DeoEst1}.

Next, we introduce the notions of symmetric and quasi-symmetric functions. For a systematic development of these, please refer to  \cite{DeoEst1} and \cite{DeoEst2}.

Let $R$ be an integral domain and $\overline{R}$ its field of fractions. Consider the polynomial ring  in $n$ variables over $R$, given by  $R\,[X]= R\,[x_1,x_2,\ldots,x_n]$. The symmetric group \sn acts in a natural way on this ring by permuting the variables.

\begin{definition}
A polynomial ${\bf f}(X) \in R\,[X]$ is said to be symmetric if it is fixed under the action of \sn. If \sn is allowed to act on $\overline{R}(X)$ in the natural way, its fixed points will be called symmetric rational functions, or simply, symmetric functions. These form a subfield $\overline{R}(X)_s$  of $\overline{R}(x)$. Furthermore, $\overline{R}(X)_s$ is generated by the $n$ elementary symmetric functions given by
\begin{eqnarray*}
  {\bf s}_{n,1}(X) &=& \sum_{i=1}^{n} x_i, \\
  {\bf s}_{n,2}(X) &=& \sum_{i<j \atop 1 \leq i,j \leq n} x_ix_j, \\
     \vdots  & &  \vdots \\
  {\bf s}_{n,n}(X) &=& x_1x_2\ldots x_n.
\end{eqnarray*}
\end{definition}

\begin{definition}
For the extension ${\Bbb F}_{q^n}/{\Bbb F}_q$, we will evaluate the elementary symmetric polynomials (resp. symmetric functions) in ${\Bbb F}_{q^n}(X)$ at $(X)=(t,\phi(t),\ldots,\phi^{n-1}(t))=(t,t^q,\ldots,t^{q^{n-1}})$. These will be called the $(n,q)$-elementary symmetric polynomials (resp. $(n,q)$-symmetric functions). For ${\bf f}(X) \in \nq\!(X)$, we will denote ${\bf f}(t,t^q,\ldots,t^{q^{n-1}})$ by $f(t)$, or, when the context is clear, by $f$.
\end{definition}

Thus the  $(n,q)$-elementary symmetric polynomials are the following:
\begin{eqnarray*}
s_{n,1}(t) &=& \sum_{0 \leq i \leq n-1} t^{q^i}, \\
s_{n,2}(t) &=& \sum_{i<j \atop 0 \leq i,j \leq n-1}t^{q^i}t^{q^j}, \\
    \vdots & & \vdots  \\
s_{n,n}(t) &=& t^{1+q+q^2+\ldots+q^{n-1}}.
\end{eqnarray*}

See \cite{DeoEst1} for a demonstration of the use of $(n,q)$-symmetric functions in splitting places of degree one in extensions of algebraic functions fields.

We now extend the notion of symmetry to get a larger class of functions that can be very effectively used to split all places of degree one in extensions of function fields. These functions are called ``$(n,q)$-quasi-symmetric.''

\begin{definition} A polynomial ${\bf f}(X)$ in $R[X]$ will be called quasi-symmetric if it is fixed by the cycle $\varepsilon = (1\; 2\;\ldots\; n) \in {\frak S}_n$.  If $\varepsilon$ is allowed to act on $\overline{R}(X)$ in the natural way, its fixed points will be called quasi-symmetric rational functions, or simply, quasi-symmetric functions. These form a subfield $\overline{R}(X)_{qs}$  of $\overline{R}(X)$.
\end{definition}

\begin{lemma} For $n>2$, there always  exist quasi-symmetric functions that are not symmetric.
\end{lemma}
{\bf Proof}.  $\langle{\varepsilon}\rangle$ has index $(n-1)!$ in $\sn$. Thus for $n >2$, the set of functions fixed by \sn is strictly smaller than those fixed by $(\varepsilon)$. For $n=2,\; \sn = (\varepsilon)$ so that the notions of symmetric and quasi-symmetric coincide.\hfill $\Box$

\begin{example} \rm \label{example:quasi-symmetric-n=3}
\rm $(n=3)$ A family of quasi-symmetric functions in three variables is given below:
$${\bf f}(x_1,x_2,x_3) = x_1{x_2}^i + x_2{x_3}^i + x_3{x_1}^i. $$
Note that for $i \neq 0 \mbox{ or }1$, these are not symmetric.
\end{example}

\begin{definition}
Consider the  extension ${\Bbb F}_{q^n}/{\Bbb F}_q$ of finite fields. We will evaluate the  quasi-symmetric polynomials (resp. quasi-symmetric functions) in ${\Bbb F}_{q^n}(X)$ at  $(X)=(t,\phi(t),\ldots,\phi^{n-1}(t))=(t,t^q,\ldots,t^{q^{n-1}})$. These will be called $(n,q)$-quasi-symmetric  polynomials (resp. $(n,q)$-quasi-symmetric functions).
\end{definition}

\begin{example}\rm
Using the three-variable quasi-symmetric functions of Example~$\ref{example:quasi-symmetric-n=3}$, we can obtain  the following $(3,q)$-quasi-symmetric functions:
$$f(t) = {\bf f}(t,t^q,t^{q^2}) = t^{1+iq} + t^{q+iq^2} + t^{q^2 + i}.$$
Again, these are not $(3,q)$-symmetric for  $i \neq 0 \mbox{ or }1$.
\end{example}

\begin{lemma} \label{lemma:qs-no-zeros}
There exist $(n,q)$-quasi-symmetric functions that have no zeros in ${\Bbb
F}_{q^n}$.
\end{lemma}
A simple method to obtain such functions is to compose irreducible polynomials over \fq, with $(n,q)$-quasi-symmetric functions.

\section{Towers where almost all rational places split completely}

In this section we construct families of towers of function fields with very good splitting behaviour. In some of the families, all rational places split completely throughout the tower, and in others, all rational places, except one, split completely throughout the tower.

\subsection{Towers of Artin-Schreier extensions}

First we begin with a tower of function fields in which all rational places split completely throughout the tower. We will denote the subfield of $F_i$ comprising $(n,q)$-quasi-symmetric functions in $x_j$ by $F_{j,qs}$ and the subfield comprising the $(n,q)$-quasi-symmetric functions of $x_j$ with \fq coefficients by $F_{j,qs}^\phi$. In particular, in $F_i$, $F_{i,qs}^\phi$ will denote the subfield of $(n,q)$-quasi-symmetric functions of $x_i$ with \fq coefficients.

\begin{theorem} \label{theorem:all-rational-split-1}
Consider the tower of function fields ${\cal F} = (F_1,F_2,\ldots)$ where $F_1 = {\Bbb F}_{q^n}(x_1)$ and for $i \geq 1$, $F_{i+1} = F_{i}(x_{i+1})$, where $x_{i+1}$ satisfies the equation
\begin{equation}
 x_{i+1}^{q^{n-1}} +  x_{i+1}^{q^{n-2}} + \ldots + x_{i+1} = \frac{g(x_i)}{h(x_i)},
\end{equation}
where $g(x_i), h(x_i) \in F_{i,qs}^\phi$, $ \frac{g(x_i)}{h(x_i)}$ is not the image of a rational function under a linear polynomial, and $h(x)$ has no zeros in \nq. Also, $deg(g(x_i)) \leq deg(h(x_i))$. Then the following hold:
\begin{enumerate}
\item All the rational places of $F_1$ split completely in all steps of the tower.
\item For every place $P$ in $T_i$ that is ramified in $T_{i+1}$, the place $P'$  in $T_{i+1}$ that lies above $P$ is unramified in $T_{i+2}$. Thus, ramification at a place cannot ``continue'' up the tower.
\end{enumerate}
\end{theorem}
{\bf Proof}. $P_\infty$ splits completely because of the condition of the degrees of $g$ and $h$. Also, the RHS is in the valuation ring at every rational place since $h$ has no zeros in \nq and $deg(g) \leq deg(h)$.  Also, its class in the residue class field is in \fq at each of these places, since the RHS is in $F_{i,qs}^\phi$. Then Proposition~\ref{proposition:Art-Sch} tells us that every rational place in $F_i$ splits completely in $F_{i+1}$. For (ii), note that if $P \in \pti$ is ramified in $T_{i+1}$, and $P'$ is a place lying above it in $T_{i+1}$, then the RHS of the equation for $x_{i+2}$ has a zero at $P'$, because of the condition on the degrees of $h$ and $g$. Thus $P'$ will be unramified in $T_{i+2}$. \hfill $\Box$

\begin{example} \rm
Consider the tower of function fields ${\cal F} = (F_1,F_2,\ldots)$ where $F_1 = {\Bbb F}_{q^3}(x_1)$, $q$ is  not a power of $2$, and for $i \geq 1$, $F_{i+1} = F_{i}(x_{i+1})$, where $x_{i+1}$ satisfies the equation
\begin{equation}
x_{i+1}^{q^2} + x_{i+1}^q + x_{i+1} = \frac{x_i^{2q^2 + 2q + 2}}{(x_i^{q^2} + x_i^q + x_i)^2 - \alpha_i},
\end{equation}
where $\alpha_i \in \fq$ is not a square.  All rational places split completely throughout the tower. Let the place $P$ of $T_1$ be a simple pole of the RHS in $T_1$ (\ie, for the case $i=1$.). Then, the place $P^{(i)}$ of $T_i$, where $i \geq 2$, which divides $P$, is a pole of $x_i$ of order $2^{i-2} \bmod q$. Also  notice that there will always exist such places, if we look at the equation over \clo. Thus the equation is absolutely irreducible at each stage.
\end{example}

\begin{theorem} \label{theorem:abelian-tower}
Consider the tower of function fields ${\cal F} = (F_1,F_2,\ldots)$ where $F_1 = {\Bbb F}_{q^n}(x_1)$, $p\neq 2$, and for $i \geq 1$, $F_{i+1} = F_{i}(x_{i+1})$, where $x_{i+1}$ satisfies the equation
\begin{equation}
 x_{i+1}^{q^{n-1}} +  x_{i+1}^{q^{n-2}} + \ldots + x_{i+1} = \frac{1}{(x_{i}^{q^{n-1}} +  x_{i}^{q^{n-2}} + \ldots + x_{i})^2 - \alpha},
\end{equation}
where $\alpha \in \fq$ is not a square. Then the following hold:
\begin{enumerate}
\item $T_i/T_1$ is an Abelian extension for $i \geq 2$.
\item All rational places split completely throughout this tower.
\item When a (non-rational) place $P \in \pti$ is ramified in $T_{i+1}$, from then on, it behaves like a rational place for splitting, and therefore splits completely further throughout the tower.
\end{enumerate}
\end{theorem}
{\bf Proof}. First we note that the equations defining the tower at each stage are indeed irreducible. For this, note that if $P$ is a place in $T_i$ that is a zero of $ (x_{i}^{q^{n-1}} +  x_{i}^{q^{n-2}} + \ldots + x_{i})^2 - \alpha)$ in $T_i$, the zero can be of degree at most two. This can be seen as follows. Let  $\sqrt{\alpha}$ be one of the square roots of $\alpha$. Then,
\begin{eqnarray*}
x_{i}^{q^{n-1}} +  x_{i}^{q^{n-2}} + \ldots + x_{i} - \sqrt{\alpha} & = &  \frac{1}{(x_{i-1}^{q^{n-1}} +  x_{i-1}^{q^{n-2}} + \ldots + x_{i-1})^2 - \alpha} - \sqrt{\alpha}, \\
& = & \frac{1 - \sqrt{\alpha}((x_{i-1}^{q^{n-1}} +  x_{i-1}^{q^{n-2}} + \ldots + x_{i-1})^2 - \alpha)}{(x_{i-1}^{q^{n-1}} +  x_{i-1}^{q^{n-2}} + \ldots + x_{i-1})^2 - \alpha}.
\end{eqnarray*}
Now note that the second derivative of the numerator of the RHS  with respect to $x_{i-1}$ is constant. The denominator is a unit at this place. Thus the zeros of the RHS can occur to at most multiplicity two. Since a similar argument holds at each stage, the valuation of the RHS at $P$ must be a power of two, which is coprime to the characteristic.  Irreducibility then follows from Proposition~\ref{proposition:Art-Sch}.
For (i), notice that the automorphisms of $T_{i+1}/T_i$ in the tower leave $x_{i+2}$ fixed, for $i \geq 1$. Further, $T_{i+1}/T_i$ is Abelian. For (ii), note that the class of the RHS in the residue field at any rational place is in \fq at any stage of the tower. And thus the defining equation splits into linear factors over the residue class field.  \hfill $\Box$

\begin{theorem} \label{theorem:wildly-ramified-all-split}
There exist wildly ramified extensions of the rational function field over non-prime fields of cardinality $> 4$ of degree equal to any power of the characteristic in which all the rational places split completely.
\end{theorem}
{\bf Proof}. For finite-separable extensions, which are not necessarily Galois, refer Theorem~\ref{theorem:all-rational-split-1}. Each extension $T_{i+1}/T_i$ has subextensions of degree equal to any arbitrary power of $p$. By an appropriate resolution of the tower, we can get the desired result.

\begin{theorem}
There exist Abelian extensions  over non-prime fields of odd characteristic of degree equal to any power of the characteristic in which all the rational places split completely.
\end{theorem}
For Abelian extensions,  Theorem~\ref{theorem:abelian-tower} says that the Galois group of the extension $T_i/T_1$ is an elementary Abelian group of exponent $p$, for $i \geq 1$. Thus, it will have normal subgroups of all indices that are powers of $p$. The result then follows by considering the fixed fields of these subgroups.

\begin{example} \rm
Consider the tower of function fields ${\cal F} = (F_1,F_2,\ldots)$ where $F_1 = {\Bbb F}_{q^3}(x_1)$, $q$ is  not a power of $2$, and for $i \geq 1$, $F_{i+1} = F_{i}(x_{i+1})$, where $x_{i+1}$ satisfies the equation
\begin{equation}
x_{i+1}^{q^2} + x_{i+1}^q + x_{i+1} = \frac{1}{(x_i^{q^2} + x_i^q + x_i)^2 - \alpha},
\end{equation}
where $\alpha \in \fq$ is not a square. In this example, all rational places split completely at all steps of the tower. Furthermore, when a (non-rational) place $P \in \pti$ is ramified in $T_{i+1}$, from then on, it behaves like a rational place for splitting, and therefore splits completely further throughout the tower.
\end{example}

\begin{theorem}
Consider the tower of function fields ${\cal F} = (F_1,F_2,\ldots)$ where $F_1 = {\Bbb F}_{q^n}(x_1)$ and for $i \geq 1$, $F_{i+1} = F_{i}(x_{i+1})$, where $x_{i+1}$ satisfies the equation
\begin{equation}
 x_{i+1}^{q^{n-1}} +  x_{i+1}^{q^{n-2}} + \ldots + x_{i+1} = \frac{g(x_i)}{h(x_i)},
\end{equation}
where $g(x_i), h(x_i) \in F_{i,qs}^\phi$, $ \frac{g(x_i)}{h(x_i)}$ is not the image of a rational function under a linear polynomial, and $h(x)$ has no zeros in \nq. Also, $deg(g(x_i)) > deg(h(x_i))$. Then all the rational places of $F_1$, except $P_\infty$, split completely in all steps of the tower.

If, in addition, we have that  $deg(g(x_i)) = deg(h(x_i)) + 1$, the the pole order of $x_i$ in the unique place lying above $P_\infty$ in $T_i$ remains one for all $i \geq 1$.
\end{theorem}

\begin{example}\rm
Consider the tower of function fields ${\cal F} = (F_1,F_2,\ldots)$ where $F_1 = {\Bbb F}_{q^3}(x_1)$, $q$ is  not a power of $2$, and for $i \geq 1$, $F_{i+1} = F_{i}(x_{i+1})$, where $x_{i+1}$ satisfies the equation
\begin{equation}
  x_{i+1}^{q^2} + x_{i+1}^q + x_{i+1} = \frac{x_i^{2q^2 + 1} + x_i^{2+q} +
x_i^{2q + q^2}}{(x_i^{q^2} + x_i^q + x_i)^2 - \alpha},
\end{equation}
where $\alpha \in \fq$ is not a square.  Here, except the unique pole $P_\infty$ of $x_1$  in $F_1$, all other rational places split completely throughout the tower. Furthermore, let $P$ be any pole of $x_2$ in $T_2$, and $P^{(n)}$ denote the unique place in $T_n$ lying above it.  Then, the  pole order of $x_n$ at $P^{(n)}$ remains constant for $n \geq 2$.
\end{example}

\begin{example} \rm
Consider the tower of function fields ${\cal F} = (F_1,F_2,\ldots)$ where $F_1 = {\Bbb F}_{q^n}(x_1)$ that is obtained as follows. $T_2 = T_1(x_1)$, where
$$ x_{2}^{q^{n-1}} +  x_{2}^{q^{n-2}} + \ldots + x_{2} = \frac{1}{(x_{1}^{q^{n-1}} +  x_{1}^{q^{n-2}} + \ldots + x_{1})^m - \alpha},
 $$
where $\alpha$ is not an $m^{th}$ power in \fq. And for $i \geq 2$, $T_{i+1} = T_i(x_{i+1})$ where $x_{i+1}$ satisfies the equation
$$ x_{i+1}^{q^{n-1}} +  x_{i+1}^{q^{n-2}} + \ldots + x_{i+1} = \frac{h(x_i)}{g(x_i)},$$
where $h(x_i),g(x_i) \in F_{i,qs}^\phi$, and $deg(h(x_i))=deg(g(x_i))+1$. Note that we are guaranteed the existence of such polynomials $h$ and $g$ by the following construction. Take any two functions $f_1$ and $f_2$ in $F_{i,qs}^\phi$ with coprime degrees $d_1$ and $d_2$ respectively (in particular, trace and norm will do). Then there exist
integers $m,n$ such that $md_1 + nd_2 = 1$. Without loss of generality, let $m$ be positive and $n$ negative. Then let $h(x_i) = i_1(f_1(x_i))$ and $g(x_i) = i_2(f_2(x_i))$, where $i_1$ and $i_2$ are irreducible polynomials over \fq of degrees $m$ and $n$ respectively. Let $P$ be any place in $T_1$ such that  $v_P({(x_{1}^{q^{n-1}} +  x_{1}^{q^{n-2}} + \ldots + x_{1})^m - \alpha}) = 1$. Then $P^{(i)}$, which is the unique place in $T_i$ dividing $P$, remains a simple pole of $x_i$ for $ i \geq 2$, ensuring irreducibility of the defining equation at each stage of the tower.
\end{example}

\subsection{Towers of Kummer extensions}

\begin{theorem} \label{theorem:tamely-ramified-all-split}
Consider the tower of function fields ${\cal F} = (F_1,F_2,\ldots)$ where $F_1 = {\Bbb F}_{q^n}(x_1)$ and for $i \geq 1$, $F_{i+1} = F_{i}(x_{i+1})$, where $x_{i+1}$ satisfies the equation
\begin{equation}
 x_{i+1}^{\frac{q^n -1}{q-1}} = \frac{g(x_i)}{h(x_i)},
\end{equation}
where $g(x_i), h(x_i) \in F_{i,qs}^\phi$, $ \frac{g(x_i)}{h(x_i)} \neq w^{\frac{q^n -1}{q-1}}, \; \forall w \in F_i$, and $g,h$ have no zeros in \nq. Also, $deg(g(x_i)) = deg(h(x_i))$. Then all the rational places split throughout the tower.
\end{theorem}
{\bf Proof}. The RHS is in the valuation ring at every rational place at $F_i,\; \forall i$ and its class in the residue class field is in $\fq \setminus \{0\}$, since $g,h$ have no zeros in \nq, and the RHS is $(n,q)$-quasi-symmetric. Then every rational place in $F_i$ splits completely in $F_{i+1}, \; \forall i$. \hfill $\Box$

\begin{example}\rm
Consider the tower of function fields ${\cal F} = (F_1,F_2,\ldots)$ where $F_1 = {\Bbb F}_{q^3}(x_1)$, $q$ is  not a power of $2$, and for $i \geq 1$, $F_{i+1} = F_{i}(x_{i+1})$, where $x_{i+1}$ satisfies the equation
\begin{equation}
  x_{i+1}^{q^2 + q + 1} = \frac{(x_i^{q^2} + x_i^q + x_i)^2 - \beta}{(x_i^{q^2} + x_i^q + x_i)^2 - \alpha},
\end{equation}
where $\alpha, \beta \in \fq$ not squares. All rational places split completely throughout the tower.
\end{example}

\begin{theorem}
Consider the tower of function fields ${\cal F} = (F_1,F_2,\ldots)$ where $F_1 = {\Bbb F}_{q^n}(x_1)$ and for $i \geq 1$, $F_{i+1} = F_{i}(x_{i+1})$, where $x_{i+1}$ satisfies the equation
\begin{equation}
 x_{i+1}^{\frac{q^n -1}{q-1}} = \frac{g(x_i)}{h(x_i)},
\end{equation}
where $g(x_i), h(x_i) \in $ are two $(n,q)$-quasi-symmetric polynomials, $ \frac{g(x_i)}{h(x_i)} \neq w^{\frac{q^n -1}{q-1}}, \forall w \in F_i$, and $g,h$ have no zeros in \nq. Also, $deg(g(x_i)) \neq deg(h(x_i))$. Then all the rational places, except possibly $P_\infty$ split throughout the tower.
\end{theorem}
{\bf Proof}. The RHS is in the valuation ring at every rational place in $F_i \; \forall i$, except possibly those dividing $P_\infty \in T_1$  and its class in the residue class field is in $\fq \setminus \{0\}$, since $g,h$ have no zeros in \nq, and the RHS is $(n,q)$-quasi-symmetric. Then every rational place in $F_i$ splits completely in $F_{i+1}, \; \forall i$. \hfill $\Box$

\begin{theorem} There exist tamely ramified extensions of arbitrarily high degree of the rational function field over any non-prime field of cardinality greater than $4$, in which all the rational places split completely.
\end{theorem}
{\bf Proof}. Consider Theorem~\ref{theorem:tamely-ramified-all-split}. We can guarantee that such $(n,q)$-quasi-symmetric functions exist, for $q>2$. Then, in the tower described in the theorem, one can go up the tower to get arbitrarily high degree extensions of the rational function field. These will not be Galois, in general. \hfill $Box$

For the towers ${\cal F}$ described in this paper in which all, or all except one, rational places split completely throughout the tower,
$\lambda({\cal F}) = 0$. This is because while the ramification in the rational places is nil, or minimal, that in the non-rational places rises quite fast, leading to a fast rise in the genus. Indeed, it seems from the known examples of towers ${\cal F}$ with $\lambda({\cal F}) > 0$ that it might be necessary to have a certain amount of ramification in the rational places, in order to have $\lambda({\cal F}) > 0$. Or at least it seems that it is not easy to control ramification in the non-rational places, and so it is better to restrict it to a few rational places alone\footnote{These statements are for towers whose first stage is a function field of genus zero.}.

\noindent{{\bf Note:}  In all the constructions given above, the property of $(n,q)$-symmetric and $(n,q)$-quasi-symmetric functions that is crucial is that they map \nq to \fq. Thus, these functions may be replaced by any other functions with this property in all the constructions. However, for such functions not to be $(n,q)$-quasi-symmetric, they must have degree atleast $q^n$ \cite{DeoEst2}.

Also, in most of the examples that appear above, we have composed the trace/norm polynomials with the irreducible polynomial $x^2 - \alpha$, where $\alpha \in \fq$ is not a square. However, we could get infinite families of further examples by using the composition $i(q(x))$, where $i(x) \in \fq[x]$ has no zeros in \fq, and $q(x)$ is a $(n,q)$-quasi-symmetric function with \fq coefficients.

\section{Towers with $\lambda({\cal F}) > 0$}

In this section, we generalize two examples of towers with $\lambda({\cal F}) > 0$ from \cite{GarSti3} to obtain two infinite families of such towers. Subfamilies attain the Drinfeld-Vladut bound.

\begin{theorem}  \label{theorem:good-family-1}
Let $q=p^n$ and $m|n, m \neq n$. Let $k_m = (p^n - 1)/(p^m -1)$.  Consider a tower of function fields in the family given by ${\cal T} = (T_1, T_2, \ldots)$, where $T_1 = \fq(x_1)$ and for $i \geq 1$, $T_{i+1} = T_i(x_{i+1})$, where $x_{i+1}$ satisfies
\begin{eqnarray*}
 x_{i+1}^{k_m} + z_i^{k_m} = b_i^{k_m}, \\
                          z_i = a_ix_i^{r_i} + b_i,
\end{eqnarray*}
where $a_i,b_i \in {\Bbb F}_{p^m} \setminus \{0\}$ for $i \geq 1$. Also $r_i$ is a power of $p, \;\forall i$. Then the following hold:
\begin{enumerate}
\item $P_\infty$ splits completely throughout the tower.
\item Every ramified place in the tower lies above a rational place in $T_1$.
\item $\lambda({\cal T}) \geq \frac{2}{q-2}$, and hence this family attains the Drinfeld-Vladut bound for $n=2, m=1$ and $q=4$.
\end{enumerate}
\end{theorem}
{\bf Proof}.
Firstly, we verify that under the hypothesis, we do indeed get a tower of function fields.   Notice that at one of the places dividing $x_1$ in $T_2$, we get a zero of $x_2$ of order not divisible by $k_m$. This implies that the RHS, for $i=1$, is not of the form $w^{k_m}$, for $w \in T_1$. Further, one of the places dividing $x_2$ in $T_3$ also exhibits the same performance, and so on up the tower. Thus, each equation is irreducible and gives us an extension.

(i) follows from the basic theory of Kummer extensions cf. \cite{Sti1}, Ch. III.7. It is important to note that linear transformations fix the place at infinity, so that it splits at each stage of the tower.

For (ii), working with residue classes, note that for ramification to take place at the $i^{th}$ step of the tower, the norm of $z_i$ should be an element of ${\Bbb F}_{p^m}$. Thus $z_i$ must be in \fq.  Since $z_i$ is obtained by a linear tranformation with \fq coefficients of a characteristic power of $x_i$, it follows that $x_i$ must be in \fq. But the relations between the variables $x_i$ and $z_{i-1}$ at the previous step of the tower then force $z_{i-1}$, and therefore $x_{i-1}$ to be in \fq. Proceeding this way to the first step of the tower, we get that $x_1 \in \fq$. Thus every ramified place in $T_i$ divides a rational place ($\neq P_\infty$) in $T_1$.

To get (iii), notice that
$$ N(T_j) > k_m^j, \mbox{ for } j \geq 1. $$
Also, the degree of the different at the \jth stage of the tower is always less than the value it would have had all $q$ finite rational places ramified from the second stage of the tower onwards. Now, using the transitivity of the different, we can say that
\begin{eqnarray*}
\degdiff(T_j/T_1) &<& q(k_m -1)[1 + k_m + \ldots + k_m^{j-2}], \\
          & < & q (k_m^{j-1} -1).
\end{eqnarray*}
Now using the Hurwitz-genus formula, it follows that
$$ g(T_j) < \frac{(q-2)(k_m^{j-1} -1)}{2}. $$
Giving us
$$\lim_{j \rightarrow \infty}  N(T_j)/g(T_j) \geq \frac{2}{q-2}. $$
\hfill $\Box$

This tower, for the case of $m=r_i=1; \, z_i = x_i+1$, first appeared in \cite{GarSti3}.

\begin{theorem} \label{theorem:good-family-2}
Let $q=p^n > 4$ and $m|n$. Let $l_m = (p^m -1)$. Consider a tower of function fields in the family given by ${\cal T} = (T_1, T_2, \ldots)$, where $T_1 = \fq(x_1)$ and for $i \geq 1$, $T_{i+1} = T_i(x_{i+1})$, where $x_{i+1}$ satisfies
\begin{eqnarray*}
x_{i+1}^{l_m} + z_i^{l_m} = 1,  \\
            z_i  = a_ix_i^{s_i} + b_i,
\end{eqnarray*}
where $a_i,b_i \in {\Bbb F}_{p^m} \setminus \{0\}$ for $i \geq 1$. Also $s_i$ is a power of $p$, $\forall i$. Then the following hold:
\begin{enumerate}
\item $P_\infty$ splits completely throughout the tower.
\item Every ramified place in the tower lies above a rational place in $T_1$ of the form $P_\gamma$, with $\gamma \in {\Bbb F}_{q^m}$.
\item $\lambda({\cal T}) \geq \frac{2}{l_m - 1}$, and hence this family attains the Drinfeld-Vladut bound for $n=2, m=1$ and $q=9$.
\end{enumerate}
\end{theorem}
{\bf Proof}. First we verify as in the proof of Theorem~\ref{theorem:good-family-1} that we do indeed get a tower of function fields. For this, note that $b_i^{l_m} = 1$. Again (i) follows from the basic theory of Kummer extensions. For (ii), we note that to have ramification at the $i^{th}$ stage of the tower, we must have that $z_i^{l_m} = 1$ implying that $z_i \in {\Bbb F}_{q^m} \setminus \{0\}$. Then by similar reasoning as in the proof of Theorem~\ref{theorem:good-family-1} , it follows that such a ramified place would divide a rational place in $T_1$ of the form $P_\gamma$, with $\gamma \in {\Bbb F}_{q^m}$. Using the Hurwitz genus formula and the transitivity property of the different along similar lines as in the proof of Theorem~\ref{theorem:good-family-1}, we get (iii). \hfill $\Box$

This tower, for the case of $s_i=m=1$, also first appeared in \cite{GarSti3}.

Following the conjecture of Elkies, it is very likely that many of  these towers are modular. In that case, there seems to be a definite relation between some modular towers and certain symmetric towers (the other optimal constructions from \cite{GarSti1} and \cite{GarSti2} are also symmetric, and are modular as shown in \cite{Elk1}). An interesting study would be to understand under what conditions can a modular tower be written down in terms of symmetric equations.

\section*{Acknowledgements}
I would like to express my deep sense of gratitude to Prof. Dennis Estes, who supervised this work,  and tragically passed away just prior to its completion. Without his constant help, I could not have made any progress whatsoever. This work is dedicated to him.

I would also like to thank Joe Wetherell for all his help in completing this work following the demise of Prof. Dennis Estes.

\end{document}